\documentclass[11pt,a4paper,leqno]{amsart}

\usepackage[T1]{fontenc}
\usepackage[english]{babel}
\usepackage{lmodern}
\usepackage{microtype}
\usepackage[margin=1in]{geometry}
\usepackage{mathtools}
\usepackage{amssymb}
\usepackage{graphicx}

\theoremstyle{plain}
\newtheorem{theorem}{Theorem}[section]
\newtheorem{lemma}[theorem]{Lemma}
\newtheorem{maintheorem}{Theorem}

\theoremstyle{definition}
\newtheorem{definition}[theorem]{Definition}

\numberwithin{equation}{section}

\newcommand{\Z}{\mathbb{Z}}
\newcommand{\F}{\mathbb{F}}
\DeclareMathOperator{\Sp}{Sp}
\DeclareMathOperator{\Mod}{Mod}
\DeclareMathOperator{\Ker}{Ker}
\let\Im\relax\DeclareMathOperator{\Im}{Im}
\DeclareMathOperator{\Int}{Int}
\DeclareMathOperator{\Arf}{Arf}
\DeclareMathOperator{\tr}{tr}

\usepackage[colorlinks=true, allcolors=blue, pdfusetitle]{hyperref}

\title[Abelianization of the genus-three Torelli commutator subgroup]{The commutator subgroup of the Torelli group of~a~genus-three surface has finitely generated~abelianization}
\author{Nikita Pshenichnyi}
\address{Lomonosov Moscow State University, Russia}
\date{}

\begin{document}

\begin{abstract}
	We prove that the commutator subgroup of the Torelli group of a genus-three surface with at most one boundary component has finitely generated abelianization. The proof combines Gaifullin's finite-generation sufficient condition with the explicit generating set for the Torelli commutator subgroup due to Brendle, Margalit, and Putman.
\end{abstract}

\maketitle
\thispagestyle{empty}

\section{Introduction and preliminaries}

Let $\Sigma_g^b$ be a compact connected oriented surface of genus $g \geqslant 3$ with $b \in \{0, 1\}$ boundary components. Its \emph{mapping class group}
\[
	\mathcal{M}_g^b \coloneqq \Mod(\Sigma_g^b)
\]
is the group of isotopy classes of orientation-preserving homeomorphisms that fix $\partial\Sigma_g^b$ pointwise, where isotopies also fix the boundary pointwise. We omit the superscript $b$ when $b=0$. The action of $\mathcal{M}_g^b$ on $H_1(\Sigma_g^b; \Z)$ preserves the algebraic intersection form and gives a surjection onto $\Sp_{2g}(\Z)$. Its kernel is the \emph{Torelli group} $\mathcal{I}_g^b$. Thus
\[
	1 \longrightarrow \mathcal{I}_g^b \longrightarrow \mathcal{M}_g^b \longrightarrow \Sp_{2g}(\Z) \longrightarrow 1.
\]

Johnson proved that $\mathcal{I}_g^b$ is finitely generated for $g \geqslant 3$ \cite{J1} and subsequently computed its~abelianization \cite{J3}; see also \cite[\S\S 1.8--1.10]{BMP}. The \emph{Johnson kernel} $\mathcal{K}_g^b$, generated by separating twists, contains $[\mathcal{I}_g^b, \mathcal{I}_g^b]$ as a subgroup of finite index \cite{J2, J3}, so the two groups are finitely generated simultaneously. Dimca and Papadima \cite[Theorem~A]{DP} proved that $H_1(\mathcal{K}_g; \mathbb{Q})$ is~finite-dimensional for $g \geqslant 4$. Ershov and He \cite[Theorem~1.1]{EH} established finite generation of $[\mathcal{I}_g^b, \mathcal{I}_g^b]$ for $g \geqslant 12$, and Church, Ershov, and Putman \cite[Theorem~A]{CEP} extended this result to~all~$g \geqslant 4$.

Gaifullin \cite[Theorem~A]{G} proved that $(\mathcal{K}_3^b)^{\mathrm{ab}}$ is finitely generated. This does not formally imply the analogous statement for $[\mathcal{I}_3^b, \mathcal{I}_3^b]^{\mathrm{ab}}$, since finite generation of abelianizations need not pass to subgroups of finite index. We prove the following.

\begin{maintheorem} \label{theorem:A}
	For $b \in \{0, 1\}$, the group $[\mathcal{I}_3^b, \mathcal{I}_3^b]^{\mathrm{ab}}$ is finitely generated.
\end{maintheorem}

The proof uses Gaifullin's finite-generation sufficient condition \cite{G}, the generating set of~Brendle, Margalit, and Putman \cite{BMP}, and Johnson's classical results \cite{J0,J1,J2,J3}.

\subsection{Separating twists and BP maps}

\begin{figure}[b]
	\centering
	\includegraphics{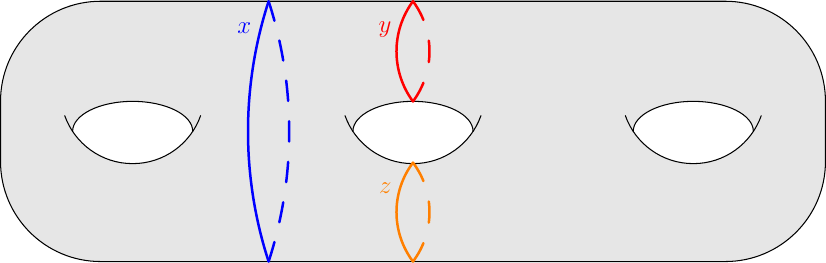}
	\caption{A separating twist $T_x$ and a BP map $T_yT_z^{-1}$.}
	\label{figure:BP}
\end{figure}

We follow the terminology of \cite[\S 1.5]{BMP}. For a simple closed curve $x$ on $\Sigma_g^b$, let $T_x$ denote the left Dehn twist about $x$. We call $x$ \emph{nontrivial} if it does not bound a disk; in particular, curves parallel to the boundary are allowed. A \emph{separating twist} is~a~Dehn twist about a nontrivial separating curve. A \emph{bounding pair \textup{(}BP\textup{)}} consists of two disjoint nonisotopic nonseparating curves $y,z$ whose union separates the surface. The corresponding \emph{BP map} is $T_yT_z^{-1}$. Separating twists and BP maps belong to $\mathcal{I}_g^b$; see Figure~\ref{figure:BP}. A mapping class is \emph{supported on} a subsurface if it has a representative equal to the identity outside that subsurface.

\subsection{Johnson homomorphisms}

Put $H \coloneqq H_1(\Sigma_g^b; \Z)$, identifying the groups for $b \in \{0, 1\}$ by~capping off the boundary. The \emph{Johnson homomorphisms} are surjections
\[
	\tau\colon \mathcal{I}_g^1 \longrightarrow U_g^1\coloneqq \Lambda^3H,\qquad
	\tau\colon \mathcal{I}_g \longrightarrow U_g \coloneqq \frac{\Lambda^3H}{\Omega\wedge H},
\]
where
\[
	\Omega \coloneqq \sum_{i = 1}^ga_i \wedge b_i \in \Lambda^2H
\]
for a symplectic basis $a_1, b_1, \ldots, a_g, b_g$ of $H$. Both target groups are free abelian. The homomorphisms are equivariant under the corresponding mapping class group actions; see \cite[Lemma~3 and \S 6]{J0}. For their construction, see \cite{J0} or \cite[\S 6.6]{FM}.

Johnson's results in \cite{J2} show that the Johnson kernel $\mathcal{K}_g^b = \Ker\tau$ is generated by separating twists (see also \cite[Theorem~A]{P}) and that the square of every separating twist belongs to $[\mathcal{I}_g^b,\mathcal{I}_g^b]$ (see also \cite[\S 1.9]{BMP}). Consequently, the kernel $\mathcal{K}_g^b / [\mathcal{I}_g^b, \mathcal{I}_g^b]$ of the induced surjection
\begin{equation} \label{eq:JohnsonAbelianization}
	\tau_{\mathrm{ab}} \colon(\mathcal{I}_g^b)^{\mathrm{ab}} \longrightarrow U_g^b
\end{equation}
is a finite elementary abelian $2$-group. Indeed, this kernel is generated by elements of order at~most two and is a subgroup of the finitely generated abelian group $(\mathcal{I}_g^b)^{\mathrm{ab}}$.

\subsection{Gaifullin's finite-generation sufficient condition} \label{section:GaifullinCondition}

For an abelian group $A$ and an element $a\in A$, write $e_a$ for the corresponding basis element of the group ring $\Z[A]$.

\begin{theorem}[{\cite[Theorem~D]{G}}]\label{theorem:Gaifullin}
	Let $A$ be a finitely generated abelian group and let $M$ be~a~finitely generated $\Z[A]$-module. Suppose that for every subgroup $B \leqslant A$ of infinite index, there exist $a_1, \ldots, a_k \in A \setminus B$, with $k \geqslant 1$, such that
	\begin{equation} \label{eq:GaifullinCondition}
		(e_{a_1} - 1) \cdots (e_{a_k} - 1)M = 0.
	\end{equation}
	Then $M$ is finitely generated as an abelian group.
\end{theorem}

We use the following equivalent formulation of \cite[Definition~1.6]{G}.

\begin{definition}
	Let $A$ be a finitely generated free abelian group and let $\mathfrak{X}$ be a nonempty subset of $A\setminus\{0\}$. Suppose that a group $\Gamma$ acts on $A$ by automorphisms. This action has the \emph{subgroup displacement property} with respect to $\mathfrak{X}$ if, for every subgroup $B\leqslant A$ of infinite index and every finite subset $X\subset\mathfrak{X}$, there exists $f\in\Gamma$ such that $f(X)\cap B=\varnothing$.
\end{definition}

Let $\kappa\colon \Lambda^3H \longrightarrow H$ be the contraction homomorphism
\[
\kappa(u\wedge v\wedge w)
\coloneqq(u\cdot v)w+(v\cdot w)u+(w\cdot u)v,
\]
where $\cdot$ denotes the algebraic intersection form. By \cite[Corollary~4.2]{G}, for $g\geqslant3$ the action of~$\Sp_{2g}(\Z)$ on $U_g^b$ has the subgroup displacement property with respect to
\begin{equation} \label{eq:DisplacementSets}
	\mathfrak{X}_g^0 \coloneqq U_g \setminus \{0\},\qquad
	\mathfrak{X}_g^1 \coloneqq \Lambda^3H \setminus \bigl((\Omega \wedge H) \cup \Ker\kappa\bigr),
\end{equation}
respectively.

\subsection{Notation and structure of the paper}

For the rest of the paper, put
\[
	D^b \coloneqq [\mathcal{I}_3^b, \mathcal{I}_3^b],\qquad
	A^b \coloneqq (\mathcal{I}_3^b)^{\mathrm{ab}},\qquad
	M^b \coloneqq (D^b)^{\mathrm{ab}} = H_1(D^b; \Z).
\]

As above, we omit the superscript when $b = 0$. We write $[f] \in A^b$ for the image of $f \in \mathcal{I}_3^b$ and $\overline{d} \in M^b$ for the image of $d \in D^b$. Our commutator convention is $[u, v] = uvu^{-1}v^{-1}$, and we use additive notation for abelian groups.

We first treat the closed surface. Define
\[
	\mathcal{Q} \coloneqq \bigl\langle\overline{T_\gamma^2} : \gamma\text{ is a nontrivial separating curve on } \Sigma_3\bigr\rangle \leqslant M.
\]
The proof of Theorem \ref{theorem:A} has three steps:

\textbf{Step \hyperref[section:Step1]{1}.} We apply Theorem~\ref{theorem:Gaifullin} to prove that $\mathcal{Q}$ is finitely generated as an abelian group.

\textbf{Step \hyperref[section:Step2]{2}.} We use \cite[Theorem~B]{BMP} to prove that $2M\subseteq\mathcal{Q}$.

\textbf{Step \hyperref[section:Step3]{3}.} We prove that $M / 2M$ is finite.

\noindent
It follows that $M / \mathcal{Q}$ is finite, being a quotient of $M / 2M$. Since $\mathcal{Q}$ is finitely generated, so is $M$.

In Section~\ref{section:Boundary}, we deduce the case $b = 1$ from the closed case using the Hochschild\,---\,Serre spectral sequence.

\section{\texorpdfstring{Finite generation of $\mathcal{Q}$}{Finite generation of Q}} \label{section:Step1}

Conjugation by $\mathcal{I}_3$ induces an action on $M$ that factors through $A$, since $D$ acts trivially on its own abelianization. Thus $M$ is a $\Z[A]$-module. It is finitely generated over $\Z[A]$: the~commutators of a finite generating set of $\mathcal{I}_3$ normally generate $D$, and their images generate $M$ as a module.

The identity $fT_\gamma^2f^{-1} = T_{f(\gamma)}^2$ shows that $\mathcal{Q}$ is invariant under conjugation by $\mathcal{M}_3$. In particular, it is a $\Z[A]$-submodule of $M$. Since $A$ is finitely generated and abelian, $\Z[A]$ is Noetherian. Hence $\mathcal{Q}$ is finitely generated over $\Z[A]$.

Choose nontrivial separating curves $x_1, \ldots, x_k$, with $k \geqslant 1$, such that $\overline{T_{x_1}^2}, \ldots, \overline{T_{x_k}^2}$ generate $\mathcal{Q}$ over $\Z[A]$. The genus-two side of each $x_i$ decomposes into a genus-one surface with two boundary components and a pair of pants, glued along curves $y_i$ and $z_i$. The BP map $C_i \coloneqq T_{y_i}T_{z_i}^{-1}$ commutes with $T_{x_i}$ and satisfies $\tau(C_i) \ne 0$; see, for instance, \cite[Corollary~3.6]{G}.

Let $B \leqslant A$ be a subgroup of infinite index. Since $\tau_{\mathrm{ab}} \colon A \to U_3$ is surjective with finite kernel, $\tau_{\mathrm{ab}}(B)$ also has infinite index. The subgroup displacement property and the equivariance of $\tau$ provide an element $f \in \mathcal{M}_3$ such that
\[
	\tau(fC_if^{-1}) \notin \tau_{\mathrm{ab}}(B),\ i = 1, \ldots, k.
\]

Here we use the surjection $\mathcal{M}_3 \to \Sp_6(\Z)$ to lift the symplectic automorphism supplied by the displacement property. Therefore, $a_i \coloneqq [fC_if^{-1}] \notin B$.

Conjugation by $f$ induces an automorphism $f_{\ast}$ of $\mathcal{Q}$ that sends $\overline{T_{\gamma}^2}$ to $\overline{T_{f(\gamma)}^2}$ for every nontrivial separating curve $\gamma$. For every $h \in \mathcal{I}_3$ and $q \in \mathcal{Q}$, this automorphism satisfies
\[
	f_{\ast}([h] \cdot q) = [fhf^{-1}] \cdot f_{\ast}(q).
\]
Therefore the classes $\overline{T_{f(x_1)}^2}, \ldots, \overline{T_{f(x_k)}^2}$ still generate $\mathcal{Q}$ over $\Z[A]$.

Since $fC_if^{-1}$ commutes with $T_{f(x_i)}$, we have
\[
	(e_{a_i} - 1)\overline{T_{f(x_i)}^2} = 0,\qquad 1 \leqslant i \leqslant k.
\]
The ring $\Z[A]$ is commutative, so
\[
	(e_{a_1} - 1) \cdots (e_{a_k} - 1)\mathcal{Q} = 0.
\]
Theorem~\ref{theorem:Gaifullin} now implies that $\mathcal{Q}$ is finitely generated as an abelian group.

\section{\texorpdfstring{The inclusion $2M \subseteq \mathcal{Q}$}{The inclusion 2M in Q}} \label{section:Step2}

We use the generating set for the Torelli commutator subgroup constructed by Brendle, Margalit, and Putman \cite{BMP}. First recall their definition of a basic commutator.

\begin{figure}[h]
	\centering
	\includegraphics{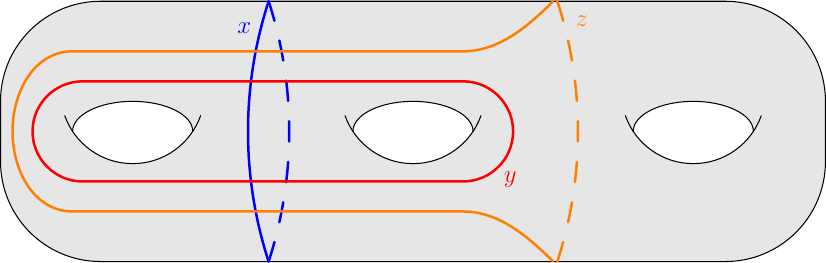}
	\caption{A basic commutator $[T_x, T_yT_z^{-1}]$.}
	\label{figure:BasicCommutator}
\end{figure}

\begin{definition}[{\cite[\S 4.1]{BMP}}]
	A commutator $[T_x, T_yT_z^{-1}]$ is \emph{basic} if $x$ is a nontrivial separating curve, $y$ and $z$ form a bounding pair, and the curves can be placed in minimal position so that
	\[
		|x \cap y| = |x \cap z| = 2
	\]
	and each of $\Sigma_g^b \setminus (x \cup y)$ and $\Sigma_g^b \setminus (x \cup z)$ has exactly two components; see Figure~\ref{figure:BasicCommutator}.
\end{definition}

\begin{theorem}[{\cite[Theorem~B]{BMP}}]\label{theorem:CommutatorGenerators}
	For $g \geqslant 3$ and $b \in \{0, 1\}$, the group $[\mathcal{I}_g^b, \mathcal{I}_g^b]$ is generated by squares of separating twists, basic commutators, and products $T_aT_bT_c$, where $a$, $b$, $c$ are disjoint nontrivial separating curves bounding a pair of pants.
\end{theorem}

The main point of this section is the following identity in $M$.

\begin{lemma} \label{lemma:BasicCommutator}
	For every basic commutator $[T_x, C]$ on $\Sigma_3$, where $C = T_yT_z^{-1}$, we have
	\begin{equation} \label{eq:BasicCommutator}
		2\overline{[T_x, C]} = \overline{T_x^2} - \overline{T_{C(x)}^2} \in \mathcal{Q}.
	\end{equation}
\end{lemma}

We first prove two auxiliary lemmas.

\begin{lemma} \label{lemma:GenusTwoSubsurface}
	For every basic commutator $[T_x, T_yT_z^{-1}]$ on $\Sigma_3$, there is a genus-two subsurface $S \subset \Sigma_3$ with one boundary component such that $x, y, z \subset \Int S$ and, after capping off $\partial S$ with a~disk, the curves $y$, $z$ become isotopic and $x$ remains nontrivial and separating.
\end{lemma}

\begin{proof}
	Place the curves as in the definition and cut $\Sigma_3$ along $y \cup z$. Both resulting components have genus one and two boundary components, and each contains two arcs of $x$. In at least one component $F$, both arcs join different boundary components. Otherwise, the arcs with endpoints on $y$ and those with endpoints on $z$ would form distinct components of $x$.

	Cutting $F$ along the first of these arcs gives a once-holed torus. The two arcs together disconnect $F$, since points on opposite sides of either arc lie in different components of $\Sigma_3 \setminus x$. Thus the second arc separates the resulting once-holed torus and cuts off a disk. It follows that $F$ contains a once-holed torus $T$, lying in its interior and disjoint from $x$, $y$, $z$.

	The torus $T$ lies on the genus-two side of $x$. Indeed, the arc of $y$ on the genus-one side is~nonseparating, since $\Sigma_3 \setminus (x \cup y)$ has exactly two components. Cutting this side along the arc gives an annulus, which cannot contain $T$.

	Set $S \coloneqq \Sigma_3 \setminus \Int T$. Capping off $\partial S$ replaces $T$ with a disk and turns $F$ into an annulus bounded by $y$, $z$, so these curves become isotopic. Both sides of $x$ then have genus one, so $x$ remains nontrivial and separating.
\end{proof}

\begin{figure}[t]
	\centering
	\includegraphics{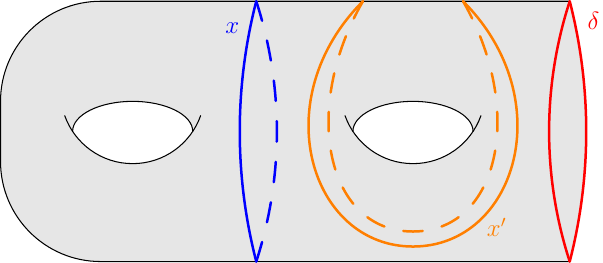}
	\caption{Curves $x$, $x^{\prime}$, $\delta$ bounding a pair of pants.}
	\label{figure:pants}
\end{figure}

\begin{lemma} \label{lemma:GeneratorsForP}
	Let $S$ be a genus-two surface with one boundary component $\delta$. Suppose that disjoint curves $x, x^{\prime} \subset \Int S$ cut $S$ into two once-holed tori and a pair of pants $\Pi$ bounded by~$x$, $x^{\prime}$, $\delta$; see Figure~\ref{figure:pants}. Let $\widehat{S} \coloneqq S \cup_{\delta} \Delta$, where $\Delta$ is a disk, and put
	\[
		\mathcal{P}_S \coloneqq \Ker\bigl(\Mod(S) \longrightarrow \Mod(\widehat S)\bigr).
	\]
	Then $\mathcal{P}_S$ has a generating set whose elements each commute with $T_x$ or $T_{x^{\prime}}$.
\end{lemma}

\begin{proof}
	The disk-pushing form of the Birman exact sequence \cite{J1} gives a central extension
	\begin{equation} \label{eq:BirmanSequence}
		1 \longrightarrow \langle T_{\delta} \rangle \longrightarrow \mathcal{P}_S \longrightarrow \pi_1(\widehat S, o) \longrightarrow 1,\ o\in\Int\Delta;\quad\text{see also \cite[\S 4.2.5]{FM}}.
	\end{equation}
	Here and below, $T_{\delta}$ denotes the twist about a curve parallel to the boundary component $\delta$.

	The subsurface $\Pi \cup_{\delta} \Delta$ is an annulus bounded by $x$ and $x^{\prime}$. Thus $\pi_1(\widehat S, o)$ is generated by four loops based at $o$: two disjoint from $x$ and two disjoint from $x^{\prime}$. Each loop has a disk-pushing lift supported away from the corresponding curve, so its lift commutes with $T_x$ or $T_{x^{\prime}}$. Together with $T_{\delta}$, which commutes with both twists, these lifts generate $\mathcal{P}_S$ by \eqref{eq:BirmanSequence}.
\end{proof}

\begin{proof}[Proof of Lemma~\ref{lemma:BasicCommutator}]
	Choose $S$ as in Lemma~\ref{lemma:GenusTwoSubsurface} and put $\delta \coloneqq \partial S$. Choose $x^{\prime}$ so that $x$ and $x^{\prime}$ cut $S$ into two once-holed tori and a pair of pants bounded by $x$, $x^{\prime}$, $\delta$. All three boundary curves of this pair of pants are nontrivial and separating in $\Sigma_3$, so Theorem~\ref{theorem:CommutatorGenerators} gives
	\[
		T_xT_{x^{\prime}}T_{\delta} \in D.
	\]

	Let $\mathcal{P}_S$ be the disk-pushing group from Lemma~\ref{lemma:GeneratorsForP}, and identify its elements with their extensions by the identity to $\Sigma_3$. Capping off $\delta$ induces an isomorphism $H_1(S; \Z) \cong H_1(\widehat{S}; \Z)$, so every element of $\mathcal{P}_S$ acts trivially on $H_1(S; \Z)$. Its extension also fixes the complementary once-holed torus. Hence $\mathcal{P}_S \subset \mathcal{I}_3$. Moreover, from Lemma \ref{lemma:GenusTwoSubsurface} it follows that $C \in \mathcal{P}_S$.

	We assert that $T_x$ fixes $\overline{[T_x, h]} \in M$ for every $h \in \mathcal{P}_S$. First take for $h$ a generator from Lemma~\ref{lemma:GeneratorsForP}. If it commutes with $T_x$, the claim is immediate. If it commutes with $T_{x^{\prime}}$, then $T_{x^{\prime}}T_{\delta}$ commutes with both $T_x$ and $h$, since $T_\delta$ is central in $\Mod(S)$. Thus $T_{x^{\prime}}T_{\delta}$ fixes $\overline{[T_x, h]}$. The product $T_xT_{x^{\prime}}T_{\delta}$ belongs to $D$ and therefore acts trivially on $M$, so $T_x$ fixes this class as well.

	Since the action of $\mathcal{I}_3$ factors through $\mathcal{I}_3^{\mathrm{ab}}$, the identity
	\[
		[T_x, gh] = [T_x, g]\,g[T_x, h]g^{-1}
	\]
	shows that this property is preserved under products and inverses. It therefore holds for every $h \in \mathcal{P}_S$, in particular for $h = C$. Finally,
	\[
		[T_x^2, C] = T_x[T_x, C]T_x^{-1}[T_x, C].
	\]
	Since $T_x$ fixes $\overline{[T_x,C]}$ and $T_x^2\in D$, we obtain
	\[
		2\overline{[T_x, C]} = \overline{[T_x^2, C]} = \overline{T_x^2} - \overline{CT_x^2C^{-1}} = \overline{T_x^2} - \overline{T_{C(x)}^2} \in \mathcal{Q}. \qedhere
	\]
\end{proof}

We now apply Theorem~\ref{theorem:CommutatorGenerators}. Classes of squares of separating twists lie in $\mathcal{Q}$ by definition, and twice the class of every basic commutator lies in $\mathcal{Q}$ by Lemma~\ref{lemma:BasicCommutator}. For a generator of the third type, the twists commute, so
\[
	2\overline{T_aT_bT_c} = \overline{T_a^2T_b^2T_c^2} = \overline{T_a^2} + \overline{T_b^2} + \overline{T_c^2} \in \mathcal{Q}.
\]
It follows that $2M \subseteq \mathcal{Q}$.

\section{\texorpdfstring{Finiteness of $M/2M$}{Finiteness of M/2M}} \label{section:Step3}

Choose a connected CW-complex $X$ with $\pi_1(X) = \mathcal{K}_3$. Since $\mathcal{K}_3 / D$ is a finite elementary abelian $2$-group, the covering corresponding to $D$ factors into a finite tower of connected double coverings. For any double covering $p\colon \widetilde{Y} \to Y$ of CW-complexes, there is a short exact sequence of cellular chain complexes
\[
	0 \longrightarrow C_{\ast}(Y; \F_2) \xlongrightarrow{\tr} C_{\ast}(\widetilde{Y}; \F_2) \xlongrightarrow{p_{\ast}} C_{\ast}(Y; \F_2) \longrightarrow 0.
\]

The transfer map $\tr\colon C_{\ast}(Y; \F_2) \to C_{\ast}(\widetilde{Y}; \F_2)$ sends each cell to the sum of its two lifts. The associated long exact sequence in homology gives
\[
	\dim H_1(\widetilde{Y}; \F_2) \leqslant 2\dim H_1(Y; \F_2).
\]

By \cite[Theorem~A]{G}, $H_1(\mathcal{K}_3; \Z)$ is finitely generated. The universal coefficient theorem therefore implies that $H_1(X; \F_2)$ is finite-dimensional. Applying the preceding estimate successively along the tower shows that $H_1(D; \F_2)\cong M / 2M$ is finite-dimensional over $\F_2$, and hence finite. Together with Sections~\ref{section:Step1} and~\ref{section:Step2}, this proves Theorem~\ref{theorem:A} for $b = 0$.

\section{The case of one boundary component} \label{section:Boundary}

Let $\delta \coloneqq \partial\Sigma_3^1$, and let $\Delta$ be the disk attached along $\delta$ to obtain $\Sigma_3$. The Birman exact sequence for the Torelli group \cite[Theorem~2.2]{P} gives
\[
	1 \longrightarrow P \longrightarrow \mathcal{I}_3^1 \longrightarrow \mathcal{I}_3 \longrightarrow 1,
\]
where $P \cong \pi_1(UT\Sigma_3)$ is the disk-pushing subgroup and $UT\Sigma_3$ is the unit tangent bundle. Restricting the capping map to commutator subgroups gives
\begin{equation} \label{eq:CommutatorExtension}
	1 \longrightarrow K \longrightarrow D^1 \longrightarrow D \longrightarrow 1,\qquad K \coloneqq P \cap D^1.
\end{equation}

\begin{lemma} \label{lemma:CappingKernel}
	The kernel in \eqref{eq:CommutatorExtension} is $K = [P, P] \langle T_\delta^2\rangle$.
\end{lemma}

\begin{proof}
	By Johnson's formula, the restriction of $\tau$ to $P$ is, up to an overall sign, the composition
\[
	P \longrightarrow\pi_1(\Sigma_3)\longrightarrow H \xlongrightarrow{\ h\,\mapsto\,\Omega \wedge h\ } \Lambda^3H,\quad\text{see \cite[Lemma~2.8]{P}}.
\]
The last map is injective, since $\kappa(\Omega \wedge h) = 2h$. As
\[
	P / \langle T_{\delta}\rangle \cong \pi_1(\Sigma_3),
\]
we obtain
\begin{equation} \label{eq:KernelInclusions}
	[P, P] \langle T_{\delta}^2\rangle \subseteq K \subseteq \Ker(\tau|_P) = [P, P]\langle T_{\delta}\rangle.
\end{equation}
The first inclusion uses $P \subset \mathcal{I}_3^1$ and $T_{\delta}^2 \in D^1$.

Choose a quadratic form $q\colon H / 2H \to \F_2$ whose associated bilinear form is the mod-$2$ intersection form and such that $\Arf(q) = 1$. Evaluating the Birman\,---\,Craggs\,---\,Johnson homomorphism at $q$ gives a homomorphism $\sigma_q\colon \mathcal{I}_3^1 \to \F_2$ satisfying $\sigma_q(T_\delta) = \Arf(q) = 1$; see \cite[\S\S 2.5--2.6]{BMP}. Since~$\sigma_q$ vanishes on $D^1$, we have $T_{\delta} \notin K$.

Thus $K$ is a proper subgroup of $[P, P] \langle T_{\delta}\rangle$. The subgroup $[P, P] \langle T_{\delta}^2\rangle$ has index at most two in this group, so the inclusions \eqref{eq:KernelInclusions} force the desired equality.
\end{proof}

Conjugation by lifts to $D^1$ defines an action of $D$ on $H_1(K; \Z)$. The action is independent of~the lifts, since inner automorphisms of $K$ act trivially on its abelianization.

\begin{lemma} \label{lemma:CoinvariantsFG}
	The coinvariant abelian group $H_1(K; \Z)_D$ is finitely generated.
\end{lemma}

\begin{proof}
	Put $N \coloneqq H_1(K; \Z)_D = K / [D^1, K]$. The groups $P$, $D^1$, and $K$ are normal in $\mathcal{M}_3^1$. Hence conjugation by $\mathcal{I}_3^1$ acts on $N$. This action factors through $A^1$, since $D^1$ acts trivially on the displayed quotient, and makes $N$ a $\Z[A^1]$-module. We will apply Theorem~\ref{theorem:Gaifullin} to this module.

	We first choose normal generators for $K$ with suitable supports. Fix $o \in \Int\Delta$ and choose simple loops $\gamma_1, \ldots, \gamma_6$ generating $\pi_1(\Sigma_3, o)$ such that any two intersect transversely only at $o$. Such a system is obtained from a dodecagon with opposite sides identified by taking the images of the segments joining the midpoints of opposite sides; see \cite[p.~182]{HM}. Take $\Delta$ sufficiently small around $o$. Choose disk-pushing lifts $p_i \in P$ supported in small regular neighborhoods of $\Delta \cup \gamma_i$, with the interior of $\Delta$ removed. Then $P$ is generated by $T_{\delta}, p_1, \ldots, p_6$. Since $T_{\delta}$ is central in $P$, Lemma~\ref{lemma:CappingKernel} shows that $T_{\delta}^2$ and $[p_i, p_j]$, $1 \leqslant i < j \leqslant 6$, normally generate $K$ in $P$.

	A regular neighborhood of $\gamma_i \cup \gamma_j$, chosen to contain $\Delta$ and the neighborhoods used for the lifts, is a once-holed torus. Removing $\Int\Delta$ gives a genus-one subsurface with two boundary components, one of which is $\delta$; this subsurface contains the supports of $p_i$ and $p_j$. The same type of subsurface supports $T_{\delta}^2$. Thus the displayed normal generators can be denoted by $k_1, \ldots, k_r$, with each $k_i$ supported on a genus-one subsurface $S_i$ with two boundary components, one of which is $\delta$. Their images generate $N$ over $\Z[A^1]$.

	The complement of each $S_i$ is a genus-two surface with one boundary component. Choose a bounding pair in this complement that bounds a genus-one subsurface, and let $C_i$ be the corresponding BP map. Then $C_i$ commutes with $k_i$, and
	\begin{equation} \label{eq:BoundaryBPImages}
		\tau(C_i) \in \mathfrak{X}_3^1 = \Lambda^3H \setminus \bigl((\Omega \wedge H) \cup \Ker\kappa\bigr).
	\end{equation}

	Indeed, Johnson's formula gives $\kappa(\tau(C_i)) = \pm c_i \ne 0$, where $c_i$ is the common homology class of~the two suitably oriented curves of the bounding pair. After capping off $\delta$, these curves still form a nontrivial bounding pair on $\Sigma_3$, so the image of $\tau(C_i)$ in $U_3$ is nonzero; see, for instance, \cite[Proposition~3.3 and Corollary~3.6]{G}.

	Let $B \leqslant A^1$ be a subgroup of infinite index. Since $\tau_{\mathrm{ab}}\colon A^1 \to \Lambda^3H$ is surjective with finite kernel, $\tau_{\mathrm{ab}}(B)$ also has infinite index. By \eqref{eq:BoundaryBPImages}, the subgroup displacement property and the equivariance of $\tau$ give an element $f \in \mathcal{M}_3^1$ such that
	\[
		\tau(fC_if^{-1}) \notin \tau_{\mathrm{ab}}(B),\qquad a_i \coloneqq [fC_if^{-1}] \notin B,\ i = 1, \ldots, r.
	\]

	Since $P$ and $K$ are normal in $\mathcal{M}_3^1$, the elements $fk_if^{-1}$ also normally generate $K$ in $P$. Their images $n_i\in N$ therefore generate $N$ over $\Z[A^1]$. Moreover, $(e_{a_i} - 1)n_i = 0$ because $fC_if^{-1}$ commutes with $fk_if^{-1}$. The ring $\Z[A^1]$ is commutative, so
	\[
		(e_{a_1} - 1) \cdots (e_{a_r} - 1)N = 0.
	\]
	Theorem~\ref{theorem:Gaifullin} implies that $N$ is finitely generated as an abelian group.
\end{proof}

\begin{proof}[Completion of the proof of Theorem~\ref{theorem:A}]
	The homological Hochschild \,---\,Serre spectral sequence of the extension \eqref{eq:CommutatorExtension} is
	\[
		E^2_{p, q} = H_p\bigl(D; H_q(K; \Z)\bigr) \Longrightarrow H_{p + q}(D^1; \Z),
	\]
	where $D$ acts on the homology of $K$ by conjugation through lifts to $D^1$. In total degree one, its second page has two terms:
	\[
		E^2_{1, 0} = H_1(D; \Z) = M,\qquad E^2_{0, 1} = H_0\bigl(D; H_1(K; \Z)\bigr) = H_1(K; \Z)_D = N.
	\]
	The first is finitely generated by the closed case, and the second by Lemma~\ref{lemma:CoinvariantsFG}.

	The differentials have bidegree $(-r, r - 1)$. No differential enters or leaves $E^2_{1, 0}$, and the only possible differential entering $E^2_{0, 1}$ is
	\[
		d^2_{2,0}\colon H_2(D; \Z) \longrightarrow N.
	\]
	No differential leaves this term. Consequently,
	\[
		E^{\infty}_{1, 0} = M,\qquad E^\infty_{0, 1} = N / \Im d^2_{2, 0}.
	\]

	The filtration of $H_1(D^1; \Z)$ supplied by the spectral sequence therefore gives a short exact sequence
	\[
		0 \longrightarrow N / \Im d^2_{2, 0} \longrightarrow H_1(D^1; \Z) \longrightarrow M \longrightarrow 0.
	\]

	Since $M$ and $N$ are finitely generated, the short exact sequence shows that $H_1(D^1; \Z) = M^1$ is finitely generated as well.
\end{proof}

\smallskip
\noindent\textsc{Acknowledgments.}~I am deeply grateful to my advisor, Alexander A. Gaifullin, for suggesting this problem and for his invaluable guidance and encouragement throughout my studies.

\smallskip
\noindent\textsc{Tool Use.}~AI was used to refine the author's English translation, proofread the final manuscript, and generate code for the figures; all AI-assisted material was reviewed by the author.

\end{document}